\documentclass{amsart}

\usepackage{amsfonts,amssymb,verbatim,amsmath,amsthm,latexsym,textcomp,amscd}
\usepackage{latexsym,amsfonts,amssymb,epsfig,verbatim}
\usepackage{amsmath,amsthm,amssymb,latexsym,graphics,textcomp}
\usepackage{graphicx}
\usepackage{color}
\usepackage{url}

\begin{document}

\newtheorem{theorem}{Theorem}[section]
\newtheorem{prop}[theorem]{Proposition}
\newtheorem{lemma}[theorem]{Lemma}
\newtheorem{cor}[theorem]{Corollary}
\newtheorem{defn}[theorem]{Definition}
\newtheorem{conj}[theorem]{Conjecture}
\newtheorem{claim}[theorem]{Claim}

\newcommand{\boundary}{\partial}
\newcommand{\C}{{\mathbb C}}
\newcommand{\D}{{\mathbb D}}
\newcommand{\cD}{\overline{\mathbb D}} 
\newcommand{\bbH}{{\mathbb H}}
\newcommand{\bbZ}{{\mathbb Z}}
\newcommand{\bbN}{{\mathbb N}}
\newcommand{\bbQ}{{\mathbb Q}}
\newcommand{\bbR}{{\mathbb R}}
\newcommand{\proj}{{\mathbb P}}
\newcommand{\lhp}{{\mathbb L}}
\newcommand{\tube}{{\mathbb T}}
\newcommand{\cusp}{{\mathbb P}}
\newcommand\AAA{{\mathcal A}}
\newcommand\BB{{\mathcal B}}
\newcommand\CC{{\mathcal C}}
\newcommand\DD{{\mathcal D}}
\newcommand\EE{{\mathcal E}}
\newcommand\FF{{\mathcal F}}
\newcommand\GG{{\mathcal G}}
\newcommand\HH{{\mathcal H}}
\newcommand\II{{\mathcal I}}
\newcommand\JJ{{\mathcal J}}
\newcommand\KK{{\mathcal K}}
\newcommand\LL{{\mathcal L}}
\newcommand\MM{{\mathcal M}}
\newcommand\NN{{\mathcal N}}
\newcommand\OO{{\mathcal O}}
\newcommand\PP{{\mathcal P}}
\newcommand\QQ{{\mathcal Q}}
\newcommand\RR{{\mathcal R}}
\newcommand\SSS{{\mathcal S}}
\newcommand\TT{{\mathcal T}}
\newcommand\UU{{\mathcal U}}
\newcommand\VV{{\mathcal V}}
\newcommand\WW{{\mathcal W}}
\newcommand\XX{{\mathcal X}}
\newcommand\YY{{\mathcal Y}}
\newcommand\ZZ{{\mathcal Z}}
\newcommand\CH{{\CC\HH}}
\newcommand\TC{{\TT\CC}}
\newcommand\EXH{{ \EE (X, \HH )}}
\newcommand\GXH{{ \GG (X, \HH )}}
\newcommand\GYH{{ \GG (Y, \HH )}}
\newcommand\PEX{{\PP\EE  (X, \HH , \GG , \LL )}}
\newcommand\MF{{\MM\FF}}
\newcommand\PMF{{\PP\kern-2pt\MM\FF}}
\newcommand\ML{{\MM\LL}}
\newcommand\PML{{\PP\kern-2pt\MM\LL}}
\newcommand\GL{{\GG\LL}}
\newcommand\Pol{{\mathcal P}}
\newcommand\half{{\textstyle{\frac12}}}
\newcommand\Half{{\frac12}}
\newcommand\Mod{\operatorname{Mod}}
\newcommand\Area{\operatorname{Area}}
\newcommand\ep{\epsilon}
\newcommand\hhat{\widehat}
\newcommand\Proj{{\mathbf P}}
\newcommand\U{{\mathbf U}}
 \newcommand\Hyp{{\mathbf H}}
\newcommand\Z{{\mathbb Z}}
\newcommand\R{{\mathbb R}}
\newcommand\Q{{\mathbb Q}}
\newcommand\E{{\mathbb E}}
\newcommand\til{\widetilde}
\newcommand\length{\operatorname{length}}
\newcommand\tr{\operatorname{tr}}
\newcommand\gesim{\succ}
\newcommand\lesim{\prec}
\newcommand\simle{\lesim}
\newcommand\simge{\gesim}
\newcommand{\simmult}{\asymp}
\newcommand{\simadd}{\mathrel{\overset{\text{\tiny $+$}}{\sim}}}
\newcommand{\ssm}{\setminus}
\newcommand{\diam}{\operatorname{diam}}
\newcommand{\pair}[1]{\langle #1\rangle}
\newcommand{\T}{{\mathbf T}}
\newcommand{\inj}{\operatorname{inj}}
\newcommand{\pleat}{\operatorname{\mathbf{pleat}}}
\newcommand{\short}{\operatorname{\mathbf{short}}}
\newcommand{\vertices}{\operatorname{vert}}
\newcommand{\collar}{\operatorname{\mathbf{collar}}}
\newcommand{\bcollar}{\operatorname{\overline{\mathbf{collar}}}}
\newcommand{\I}{{\mathbf I}}
\newcommand{\tprec}{\prec_t}
\newcommand{\fprec}{\prec_f}
\newcommand{\bprec}{\prec_b}
\newcommand{\pprec}{\prec_p}
\newcommand{\ppreceq}{\preceq_p}
\newcommand{\sprec}{\prec_s}
\newcommand{\cpreceq}{\preceq_c}
\newcommand{\cprec}{\prec_c}
\newcommand{\topprec}{\prec_{\rm top}}
\newcommand{\Topprec}{\prec_{\rm TOP}}
\newcommand{\fsub}{\mathrel{\scriptstyle\searrow}}
\newcommand{\bsub}{\mathrel{\scriptstyle\swarrow}}
\newcommand{\fsubd}{\mathrel{{\scriptstyle\searrow}\kern-1ex^d\kern0.5ex}}
\newcommand{\bsubd}{\mathrel{{\scriptstyle\swarrow}\kern-1.6ex^d\kern0.8ex}}
\newcommand{\fsubeq}{\mathrel{\raise-.7ex\hbox{$\overset{\searrow}{=}$}}}
\newcommand{\bsubeq}{\mathrel{\raise-.7ex\hbox{$\overset{\swarrow}{=}$}}}
\newcommand{\tw}{\operatorname{tw}}
\newcommand{\base}{\operatorname{base}}
\newcommand{\trans}{\operatorname{trans}}
\newcommand{\rest}{|_}
\newcommand{\bbar}{\overline}
\newcommand{\UML}{\operatorname{\UU\MM\LL}}
\newcommand{\EL}{\mathcal{EL}}
\newcommand{\tsum}{\sideset{}{'}\sum}
\newcommand{\tsh}[1]{\left\{\kern-.9ex\left\{#1\right\}\kern-.9ex\right\}}
\newcommand{\Tsh}[2]{\tsh{#2}_{#1}}
\newcommand{\qeq}{\mathrel{\approx}}
\newcommand{\Qeq}[1]{\mathrel{\approx_{#1}}}
\newcommand{\qle}{\lesssim}
\newcommand{\Qle}[1]{\mathrel{\lesssim_{#1}}}
\newcommand{\simp}{\operatorname{simp}}
\newcommand{\vsucc}{\operatorname{succ}}
\newcommand{\vpred}{\operatorname{pred}}
\newcommand\fhalf[1]{\overrightarrow {#1}}
\newcommand\bhalf[1]{\overleftarrow {#1}}
\newcommand\sleft{_{\text{left}}}
\newcommand\sright{_{\text{right}}}
\newcommand\sbtop{_{\text{top}}}
\newcommand\sbot{_{\text{bot}}}
\newcommand\sll{_{\mathbf l}}
\newcommand\srr{_{\mathbf r}}
\newcommand\geod{\operatorname{\mathbf g}}
\newcommand\mtorus[1]{\boundary U(#1)}
\newcommand\A{\mathbf A}
\newcommand\Aleft[1]{\A\sleft(#1)}
\newcommand\Aright[1]{\A\sright(#1)}
\newcommand\Atop[1]{\A\sbtop(#1)}
\newcommand\Abot[1]{\A\sbot(#1)}
\newcommand\boundvert{{\boundary_{||}}}
\newcommand\storus[1]{U(#1)}
\newcommand\Momega{\omega_M}
\newcommand\nomega{\omega_\nu}
\newcommand\twist{\operatorname{tw}}
\newcommand\modl{M_\nu}
\newcommand\MT{{\mathbb T}}
\newcommand\Teich{{\mathcal T}}
\newcommand{\df}{\hbox{Diff}^{\omega}(S^1)}
\newcommand{\dfal}{\hbox{Diff}^{\omega}_{\alpha}(S^1)}
\newcommand{\dinf}{{\mathbb D}_{\infty}}
\renewcommand{\Re}{\operatorname{Re}}
\renewcommand{\Im}{\operatorname{Im}}

\title{Loewner evolution of hedgehogs and 2-conformal measures of circle maps}

\author{Kingshook Biswas }

\date{}


\begin{abstract} Let $f$ be a germ of holomorphic
diffeomorphism with an irrationally indifferent fixed point at the
origin in $\C$ (i.e. $f(0) = 0, f'(0) = e^{2\pi i \alpha},
\alpha \in \R - \Q$). Perez-Marco (Fixed points and circle maps, Acta Math. {\bf 179} (2) 1997, 243-294)
showed the existence of a
unique continuous monotone one-parameter family of nontrivial invariant full continua containing the
fixed point called Siegel compacta, and gave a correspondence between germs and families $(g_t)$ of circle
maps obtained by conformally mapping the complement of these compacts to the
complement of the unit disk. The family of circle maps $(g_t)$ is the orbit
of a locally-defined semigroup $(\Phi_t)$ on the space of analytic circle maps which we
show has a well-defined infinitesimal generator $X$. The explicit form of $X$ is obtained
by using the Loewner equation associated to the family of hulls $(K_t)$. We show that the Loewner measures
$(\mu_t)$ driving the equation are 2-conformal measures on the circle for the circle maps $(z \mapsto \overline{g_t(\overline{z})})$.

\smallskip

\begin{center}

{\em AMS Subject Classification: 37F50}

\end{center}

\end{abstract}

\maketitle

\overfullrule=0pt

\tableofcontents

\section{Introduction}

\bigskip

A germ $f(z) = e^{2\pi i \alpha}z + O(z^2), \alpha \in \mathbb{R -
Q}/\mathbb{Z}$ of holomorphic diffeomorphism fixing the origin in
$\mathbb{C}$ is said to be linearizable if it is analytically
conjugate to the rigid rotation $R_{\alpha}(z) = e^{2\pi i \alpha}
z$. The number $\alpha$ is called the rotation number of $f$, and
the maximal domain of linearization is called the Siegel disk of
$f$. The linearizability of $f$ is dependent on the arithmetic of
$\alpha$, and the optimal arithmetic condition for linearizability
in this setting is given by the well-known Brjuno condition (see \cite{siegel}, \cite{brjuno},
\cite{yoccoz}).

\medskip

Perez-Marco proved the existence of a unique, strictly increasing Hausdorff continuous
family $(K_t)_{t > 0}$ of non-trivial invariant full continua containing the
fixed point called {\it Siegel compacts} (\cite{perezmcirclemaps}), where $K_t \to \{0\}$
as $t \to 0$, and $K_t$ can be described as the connected component containing the origin
of the set of non-escaping points in the closed disc ${\cD}_t$ of radius $t$ around the
origin. When $f$ is
non-linearizable these are called {\it hedgehogs}. The topology and dynamics of hedgehogs
has been studied by Perez-Marco (\cite{perezmtopology},
\cite{perezmdynamics}), who also developed techniques using  "tube-log Riemann
surfaces" (\cite{bpm1}, \cite{bpm2}, \cite{bpm3})
for the construction of interesting examples \cite{perezmens}, \cite{perezminvent},
\cite{perezmsmooth} of indifferent germs and hedgehogs, which were also used by the author
(\cite{biswas1}, \cite{biswas2}, \cite{biswas3})
and Cheritat (\cite{cheritat}) to construct further examples.

\medskip

\noindent{\bf Notation:} Throughout $\dinf$ will denote the complement in $\hat{\C}$ of the closed unit disk, $\dinf = \hat{\C} - \cD$,
$r : S^1 \to S^1$ will denote complex conjugation, $r(\xi) = \overline{\xi}$, and $\frac{d^+}{dt}$ will denote right-hand
derivative.

\medskip

The construction of Perez-Marco from \cite{perezmcirclemaps} associates to a pair
$(f, K)$ an analytic circle diffeomorphism $g$, where $f$ is a germ with an irrationally indifferent
fixed point at the origin and $K$ is a Siegel compact of $f$, by considering a conformal map
$\psi$ from the complement of $K$ to the complement of the closed unit disk,
$\psi : \hat{\C} - K \to \dinf$, such that $\psi(\infty) = \infty$. Conjugating $f$ by $\psi$
gives a holomorphic diffeomorphism $g = \psi \circ f \circ \psi^{-1}$ in an annulus in $\dinf$
having $S^1$ as a boundary component, and $g$ is shown to extend across $S^1$ to an analytic
circle diffeomorphism defined in a neighbourhood of $S^1$ such that the rotation numbers of $g$ and $f$
are equal, $\rho(g) = \rho(f) = \alpha$.

\medskip

Invariant compacts for $g$ containing $S^1$ then correspond to invariant compacts for
$f$ containing $K$, and the theorem on existence and uniqueness of Siegel compacts then gives
the existence and uniqueness of a continuous strictly increasing $1$-parameter family of {\it Herman compacts} $(A_t)_{0 \leq t < \epsilon}$
for the circle map $g$ with $A_0 = S^1$, where a Herman compact for $g$ is a connected totally invariant compact
$A$ containing $S^1$ such that $\hat{\C} - A$ has two components (such a compact is necessarily invariant
under reflection in the unit circle as well). The construction which produces $g$ from the pair
$(f, K)$ can similarly be applied to each pair $(g, A_t)$ to give a one-parameter family of analytic
circle diffeomorphisms $(g_t)_{t \geq 0}$ with $g_0 = g$, by conjugating $g$ by a conformal map
$\psi_t : \hat{\C} - (\cD \cup A_t) \to \dinf$, normalized so that $\psi_t(\infty) = \infty,
\psi'_t(\infty) > 0$. As before the map $g_t = \psi_t \circ g \circ \psi^{-1}_t$
defined in a one-sided neighbourhood of $S^1$ extends across to give an analytic circle diffeomorphism
$g_t$ with rotation number equal to that of $g$.

\medskip

The family $A_t$ being Hausdorff continuous and strictly increasing, it is possible to reparametrize the family
so that the compacts $\cD \cup A_t$ have logarithmic capacity $t$, meaning that the expansion of $\psi_t$ near
$z = \infty$ is of the form
$$
\psi_t(z) = e^{-t} z + a_0(t) + \frac{a_1(t)}{z} + \dots
$$
Let $\df$ and $\dfal$ denote the space of analytic circle diffeomorphisms and the
subspace of those diffeomorphisms with fixed rotation number $\alpha \in (\R - \Q)/\Z$ respectively.
We obtain for $t \geq 0$ a locally defined family of maps $\Phi_t : \DD_t \subset \dfal \to \dfal$ mapping
a circle map $g$ to the circle map $g_t$ obtained from the above construction.
Here the domain $\DD_t$ of $\Phi_t$ consists of those circle maps in $\dfal$ which have a Herman compact $A$
such that $\cD \cup A$ has logarithmic capacity $t$. It is not hard to show that in fact the maps $\Phi_t$
form a continuous semigroup, namely $\Phi_0 = id$,
$$
\Phi_s \circ \Phi_t = \Phi_{s+t}
$$
whenever the composition is defined, and the orbits $t \mapsto \Phi_t(g)$ give continuous curves in $\dfal$.

\medskip

We show in fact that the semigroup $(\Phi_t)_{t \geq 0}$
has an infinitesimal generator $X$, meaning the curves $t \mapsto \Phi_t(g)$ are right-hand differentiable
in $t$ and
$$
\frac{d^+}{dt}\Phi_t(g) = X(\Phi_t(g))
$$
Since the space $\dfal$ does not carry any obvious differentiable structure, the sense in which
these assertions hold is made precise in the statements of the theorems below.

\medskip

The form of the infinitesimal generator $X$ is obtained by studying the {\it Loewner equation} associated
to the family of hulls $\cD \cup A_t$ (we recall in section 2 the basic facts about the Loewner equation
which we will be needing). The maps $\phi_t := \psi^{-1}_t : \dinf \to \hat{\C} - (\cD \cup A_t)$ form a
{\it Loewner chain}, and it is known that $t \mapsto \phi_t(z,t)$ is absolutely continuous for each
fixed $z \in \dinf$. Moreover, for almost every $t$, the right-hand derivatives
$$
\chi_t(z) = \frac{d^+}{ds}_{|s = t} (\phi^{-1}_t \circ \phi_s)(z)
$$
exist and the functions $H_t(z) = \chi_t(z)/z$ are
given by the Herglotz transforms on $\dinf$ of a family of probability measures
on the unit circle $(\mu_t)$,
called the {\it driving} or {\it Loewner measures} of the Loewner equation.
Here by the Herglotz transform on $\dinf$ of a
probability measure
$\mu$ on $S^1$ we mean the holomorphic function $H = \mathcal{H}\mu$ in $\dinf$ with positive
real part and satisfying $H(\infty) = 1$, defined by
$$
(\HH \mu)(z) = \int_{S^1} \frac{\xi + 1/z}{\xi - 1/z} d\mu(\xi)
$$
(the classical Herglotz Theorem asserts that any holomorphic function $H$ on $\dinf$
satisfying $\Re H > 0, H(\infty) = 1$ is of this form for a unique probability measure $\mu$ on $S^1$).
Douady and Yoccoz have shown in \cite{douadyyoccoz}
the existence and uniqueness of an $s$-conformal measure $\mu = \mu_{s,g}$ for any $C^2$ circle diffeomorphism
$g$ with irrational rotation number and any $s \in \R$. By an $s$-conformal
measure for a circle diffeomorphism $g$
we mean a probability measure $\mu$ on $S^1$ such that
$$
\mu(g(E)) = \int_E |g'(x)|^s d\mu(x)
$$
for all measurable sets $E \subset S^1$.

\medskip

We show that for any $g$ in $\dfal$, the associated Loewner chain $(\phi_t)$ is differentiable
for {\bf all} $t$ uniformly
for $z$ in any compact subset of $\dinf$, and the associated Loewner measures are
the pull-backs by complex conjugation of the $2$-conformal measures of the circle diffeomorphisms $g_t = \Phi_t(g)$:

\medskip

\begin{theorem} \label{meastwoconform} (Loewner measures are $2$-conformal measures). For any $g$ in $\dfal$
and any $t \geq 0$, the right-hand derivative of the Loewner chain $(\phi_t)$
exists (uniformly for $z$ in any compact subset of $\dinf$) and is given by
$$
\chi_t(z) = \frac{d^+}{ds}_{|s = t} (\phi^{-1}_t \circ \phi_s)(z) = z \cdot (\HH r^* \mu_{2, g_t})(z)
$$
where $\mu_{2, g_t}$ is the unique $2$-conformal measure of $g_t$.
\end{theorem}

\medskip

The existence and form of the infinitesimal generator $X$ of the semigroup $(\Phi_t)$ is stated as follows:

\medskip

\begin{theorem} \label{infgenexists} (Infinitesimal generator of the semigroup).
For any $g \in \dfal$, there exists a function $X(g)$ holomorphic in a
neighbourhood $V$ of $S^1$ such that for $t > 0$ small the circle maps $\Phi_t(g)$ are defined
in $V$ and
$$
\lim_{t \to 0^+} \frac{\Phi_t(g) - g}{t} = X(g)
$$
uniformly on compacts in $V$. The holomorphic function $X(g)$ is given (in $V \cap \dinf$) by
$$
X(g) = g' \cdot \chi - \chi \circ g
$$
where $\chi(z) = z \cdot (\HH r^*\mu_{2,g})(z)$.
\end{theorem}

\medskip

We can think of the curves $(t \mapsto \Phi_t(g))$ as integral curves of a vector field $X$
on the space $\dfal$. The next theorem asserts the uniqueness in forward time for these integral
curves:

\medskip

\begin{theorem} \label{fwduniq} (Uniqueness in forward time of integral curves). If $(g_t)_{0 \leq t < \epsilon} \subset \dfal$
is a family of circle diffeomorphisms holomorphic in a neighbourhood $V$ of $S^1$ such that
the right-hand derivatives exist and satisfy
$$
\frac{d^+}{ds}_{|s = t} g_s(z) = X(g_t)(z)
$$
uniformly on compacts in $V$, then $g_t = \Phi_t(g_0)$ for $0 \leq t < \epsilon$.
\end{theorem}

\medskip

\begin{defn} [Germ of integral curve] An integral curve (in the sense of the
previous theorem) $(g_t)_{-\infty < t < t_0} \subset \dfal$
defined for all times $t$ less than some $t_0$ is said to be a {\it backward integral curve} of $X$.
We say that two backward integral curves $(g^1_t)_{-\infty < t < t_1}, (g^2_t)_{-\infty < t < t_2}$ of $X$
define the same germ of integral curve of $X$ if $g^1_t = g^2_t$ for all $t < t_0$ for some
$t_0 < t_1, t_2$.
\end{defn}

\medskip

It follows from the above that if $(g^1_t)_{-\infty < t < t_0}, (g^2_t)_{-\infty < t < t_0}$
are two backward integral curves defining the same germ of integral curve
then in fact $g^1_t = g^2_t$ for {\bf all} $t < t_0$.
Any germ $f$ with an irrationally indifferent fixed point gives a family of Siegel
compacts $(K_t)_{-\infty < t < t_0}$ (parametrizing the compacts $K_t$ by their logarithmic
capacities), and hence (by applying the germs to circle maps construction to the pairs $(f, K_t)$)
gives a family of circle maps $(g_t = g^f_t)_{-\infty < t < t_0}$ which it is easy to see is
a backward integral curve. We show that conversely any backward integral curve arises in this
way from a germ $f$. Denoting the space of germs with rotation number $\alpha$ by
$\hbox{Diff}_{\alpha}(\C, 0)$ we have:

\medskip

\begin{theorem} \label{germintegral} (Germs of diffeomorphisms and germs of integral curves).
For any backward integral curve $(g_t)_{-\infty < t < t_0} \subset \dfal$ of $X$,
we have $g_t \to R_{\alpha}$ uniformly on $S^1$ as $t \to -\infty$.

\smallskip

\noindent The map $f \mapsto [(g^f_t)]$ gives a one-to-one
correspondence between $\hbox{Diff}_{\alpha}(\C, 0)$ and germs of integral curves of $X$.
\end{theorem}

\medskip

Finally in section 8 we describe some further results that can be obtained in the case of analytically
linearizable circle maps and germs.

\medskip

\noindent{\bf Acknowledgements} The author thanks Ricardo
Perez-Marco for helpful discussions.

\bigskip

\section{Boundary values of the Herglotz transform}

\medskip

The Herglotz transform of a probability measure $\mu$ on $S^1$ is the holomorphic function $\HH \mu$ in $\D$ defined by
$$
(\HH \mu)(w) = \int_{S^1} \frac{\xi + w}{\xi - w} d\mu(\xi)
$$
The real and imaginary parts of the Herglotz transform $\HH \mu = \PP \mu + i \QQ \mu$ are given by the Poisson
 and conjugate Poisson transforms,
\begin{align*}
(\PP \mu)(w) & = \int_{S^1} \Re \frac{\xi + w}{\xi - w} d\mu(w) = \int_{S^1} \frac{1 - |w|^2}{|\xi - w|^2} d\mu(\xi) \\
(\QQ \mu)(w) & = \int_{S^1} \Im \frac{\xi + w}{\xi - w} d\mu(w) = \int_{S^1} \frac{2 \Im \overline{\xi}w}{|\xi - w|^2} d\mu(\xi) \\
\end{align*}
The radial limits of these harmonic functions exist for a.e. $\xi \in S^1$ with respect to Lebesgue measure $\lambda$,
and the decomposition of the measure $\mu$ into absolutely continuous and singular parts with respect to Lebesgue
measure $\mu = f d\lambda + \mu_s$ can be recovered from these radial limits as follows:

\medskip

\begin{theorem} \label{hreal} (Fatou). For Lebesgue-a.e. $\xi \in S^1$. the radial limit of $(\PP\mu)(w)$ exists and equals $f(\xi)$.
\end{theorem}

\medskip

\begin{theorem} \label{himaginary} (Poltoratski \cite{poltoratski}). Let $Q$ denote the radial limit of $\QQ \mu$ on $S^1$.
As $t \to +\infty$, the measures $\frac{\pi}{2}t 1_{\{|Q| > t\}} d\lambda$
converge weakly to $\mu_s$.
\end{theorem}

\medskip

Since we will be dealing with functions defined in $\dinf$, we will refer to
the holomorphic function in $\dinf$ defined by $H(z) = \int_{S^1} \frac{\xi + 1/z}{\xi - 1/z} d\mu(\xi)$
as the Herglotz transform of $\mu$ on $\dinf$. Then $H(z) = \tilde{H}(w)$, where $w = 1/z$ and $\tilde{H} = \HH \mu$
is the Herglotz transform defined above which is holomorphic in $\D$.
Appropriate versions of the above theorems hold for the boundary
values of the function $H$. Let the boundary
values of $H, \tilde{H}$ be $P+iQ, \tilde{P}+i\tilde{Q}$ respectively. Then $P = \tilde{P} \circ r,
Q = \tilde{Q} \circ r$, so for $\mu = f d\lambda + \mu_s$ and $\xi \in S^1$ we have
$$
f \circ r = \tilde{P} \circ r = P
$$
and
\begin{align*}
r^* \mu_s & = \lim_{t \to +\infty} r^*\left(\frac{\pi}{2} t 1_{\{|\tilde{Q}| > t\}} d\lambda\right) \\
          & = \lim_{t \to +\infty} \left(\frac{\pi}{2} t 1_{\{|\tilde{Q} \circ r| > t\}} d(r^* \lambda)\right) \\
          & = \lim_{t \to +\infty} \left(\frac{\pi}{2} t 1_{\{|Q| > t\}} d\lambda \right) \\
\end{align*}
weakly.

\bigskip

\section{The Loewner equation}

\medskip

A {\it hull} is a connected, full, non-trivial compact $K$ in $\C$ containing
the origin. Its complement $\Omega$ in $\hat{\C}$ is then a simply connected
domain containing $\infty$, so there is a conformal map from the complement
of the closed unit disk, $\phi : \dinf \to \Omega = \hat{\C} - K$ such that
$\phi(\infty) = \infty$, which is unique when normalized to satisfy $\phi'(\infty) > 0$. The map
$\phi$ has an expansion around $z = \infty$ of the form
$$
\phi(z) = e^{c(K)} z + a_0 + \frac{a_1}{z} + \dots
$$
where the real number $c(K)$ is called the {\it logarithmic capacity} of the hull $K$
(note that the closed disk of radius $R$ then has logarithmic capacity $\log R$).

\medskip

Given a strictly increasing family of hulls $(K_t)_{-\infty < t \leq t_0}$, if the domains
$\Omega_t = \hat{\C} - K_t$ are continuous for the Caratheodory topology (which holds if the
family of hulls is Hausdorff continuous for example) and $K_t \to \{0\}$ as $t \to -\infty$,
then one can continuously reparametrize the family by logarithmic capacity so
that the associated conformal maps $\phi_t : \dinf \to \Omega_t$ satisfy
$$
\phi_t(z) = e^{t} z + a_0(t) + \frac{a_1(t)}{z} + \dots
$$
The family of conformal maps $(\phi_t)$ is called a {\it Loewner chain}.

\medskip

For each compact $K \subset \dinf$,
for each $z$ in $K$ the map $t \mapsto \phi_t(z)$ is $C_K$-Lipschitz in $t$ for some constant
$C_K$ only depending on $K$ (\cite{pomm}, \S 6.1).  Thus for each $z$ the derivative
of $\phi_t(z)$ with respect to $t$ exists for all $t$ outside some null set $E_z$, but in fact one
can choose a null set $E$ independent of $z$ such that the derivative with respect to $t$ exists for
all $z$ in $\dinf$. Moreover this derivative is of the form
$$
\frac{d^+}{ds}_{|s = t} (\phi^{-1}_t \circ \phi_s)(z) = zH_t(z)
$$
where $H_t$ is a holomorphic function on $\dinf$ satisfying $\Re H_t > 0, H_t(\infty) = 1$,
and hence can be written as the Herglotz transform on $\dinf$ of a
unique probability measure $\mu_t$ on $S^1$,
$$
H_t(z) = \int_{S^1} \frac{\xi + 1/z}{\xi - 1/z} d\mu_t(\xi)
$$
The measures $(\mu_t)$ are called the driving or Loewner measures for the Loewner chain $(\phi_t)$.

\bigskip

\section{Conformal measures of analytic circle diffeomorphisms}

\medskip

A conformal measure of dimension $\delta$ for a holomorphic map $f$ is a
finite measure $\mu$ on a $f$-invariant compact $K$ such that
$$
\mu(f(E)) = \int_E |f'(z)|^{\delta} d\mu(z)
$$
for all measurable sets $E \subset K$, or, when $f$ is a diffeomorphism,
$$
f^* \mu = |f'|^{\delta} d\mu
$$
where $f^* \mu$ is defined by $(f^* \mu)(E) = \mu(f(E))$. Conformal measures have been studied extensively for
rational maps and Kleinian groups (as Patterson-Sullivan measures), where existence is usually proved by means of transfer
operators for hyperbolic dynamical systems. While the machinery of transfer operators is
 not available for circle diffeomorphisms, Douady and Yoccoz prove (\cite{douadyyoccoz})
 nonetheless the existence and uniqueness of $s$-conformal measures $\mu_{s,g}$
for any $C^2$ circle diffeomorphism $g$ with irrational rotation number for every $s \in \R$.
Moreover the measures $\mu_{s,g}$ depend continuously on $s, g$ for the weak topology
on measures and the $C^1$ topology on circle diffeomorphisms.

\bigskip

\section{The semigroup $(\Phi_t)_{t \geq 0}$ on $\dfal$}

\medskip

Given an analytic circle diffeomorphism $g$ with irrational rotation number $\alpha$,
for $\epsilon > 0$ small, there is a unique, strictly increasing, Hausdorff continuous
family of Herman compacts
$(A_t)_{0 \leq t < \epsilon}$ totally invariant under $g$ such that the hulls $\cD \cup A_t$
have logarithmic capacity $t$. Let $(\Omega_t)_{0 \leq t < \epsilon}$ be the
family of decreasing simply connected domains $\Omega_t = \hat{\C} - (\cD \cup A_t)$, and
let $\phi_t : \dinf \to \Omega_t$ be conformal maps normalized so that $\phi_t(\infty) = \infty,
\phi'_t(\infty) = e^t$. We let $\psi_t = \phi^{-1}_t : \Omega_t \to \dinf$.

\medskip

The map $g_t := \psi_t \circ g \circ \psi^{-1}_t$ is holomorphic in an annulus contained in
$\dinf$ with one boundary component equal to $S^1$, and (as shown in \cite{perezmcirclemaps})
extends analytically across $S^1$ to give an analytic circle diffeomorphism with rotation number equal
to $\alpha$. This defines a map
\begin{align*}
\Phi_t : \DD_t \subset \dfal & \to \dfal \\
                           g & \mapsto g_t \\
\end{align*}
where $\DD_t$ is the set of $g$ in $\dfal$ having a Herman compact $A$ such that $\cD \cup A$ has
logarithmic capacity $t$.

\medskip

For $g$ in $\dfal$, and $t \geq 0$, let $A_t(g)$ be a Herman compact for $g$ such that $\cD \cup A_t(g)$
has logarithmic capacity $t$ (if such a Herman compact exists), let $\Omega_t(g) := \hat{\C} - (\cD \cup A_t(g))$,
and let $\psi_t(g) : \Omega_t(g) \to \dinf$ be the corresponding normalized conformal map.

\medskip

\begin{prop} For $s,t \geq 0$, $\Phi_t(\DD_{s+t}) \subset \DD_s$, and $\Phi_s \circ \Phi_t = \Phi_{s+t}$
on $\DD_{t+s}$.
\end{prop}

\medskip

\noindent{\bf Proof:} For $g \in \DD_{s+t}$ where $s, t \geq 0$, $\psi_{t}(g)$ maps $\Omega_{s+t}(g)$ conformally
to a domain $\Omega_1$ such that the compact $A_1$ given by the complement of $\Omega_1$
and its reflection in the unit circle is invariant under $\Phi_t(g)$,
 hence $\Omega_1 = \Omega_{s'}(\Phi_t(g))$ for some $s' \geq 0$. The map $\psi_{s'}(\Phi_t(g))$
 maps $\Omega_1$ conformally to $\dinf$, so the composition $\psi := \psi_{s'}(\Phi_t(g)) \circ \psi_t(g)$
 maps $\Omega_{s+t}(g)$
conformally to $\dinf$ and satisfies the normalizations $\psi(\infty) = \infty,
\psi'(\infty) = e^{s'} e^t = e^{s'+t} > 0$. It follows from uniqueness of the normalized
conformal mapping that $\psi = \psi_{s+t}(g)$, so $e^{s'+t} = e^{s+t}$ and hence $s' = s$.
Thus $\psi_{s}(\Phi_t(g)) \circ \psi_t(g) = \psi_{s+t}(g)$, from which it follows that
$(\Phi_s \circ \Phi_t)(g) = \Phi_{s+t}(g)$. $\diamond$

\bigskip

\section{Infinitesimal generator of the semigroup}

\medskip

Let $g$ be an analytic circle diffeomorphism with irrational rotation number $\alpha$.
Let $(A_t)_{0 \leq t < \epsilon}$ be the family of Herman compacts of $g$, and $(\HH_t = \cD \cup A_t)_{0 \leq t < \epsilon}$
the associated family of hulls. Let $\Omega_t = \hat{\C} - \HH_t$, and let
$\phi_t : \dinf \to \Omega_t, \psi_t = \phi^{-1}_t$ be normalized conformal mappings. Then
the functions $\{ (\phi_t - id)/t : 0 \leq t < \epsilon \}$ form a normal family (\cite{pomm}, \S 6.1).

\bigskip

\begin{lemma} \label{psidot} Let $\{ t_n \}$ be a sequence decreasing monotonically to $0$ such that
$(\phi_{t_n} - id)/t_n$ converges uniformly on compacts of $\dinf$
to a holomorphic function $\chi$. Then $(\psi_{t_n} - id)/t_n$ converges uniformly on compacts
in $\dinf$ to $(-\chi)$.
\end{lemma}

\medskip

\noindent{\bf Proof:} Fix a compact $K \subset \dinf$, then for $n$ large $K \subset \Omega_{t_n}$
and the maps $\psi_{t_n}$ are defined on $K$. We have $\phi_{t_n} = id + t_n \chi + o(t_n)$
uniformly on $K$. Then $id = \phi_{t_n} \circ \psi_{t_n}$ gives
$$
id = (id + t_n \chi + o(t_n)) \circ \psi_{t_n} = \psi_{t_n} + t_n \chi \circ \psi_{t_n} + o(t_n),
$$
so
$$
\frac{\psi_{t_n} - id}{t_n} = -\chi \circ \psi_{t_n} + o(1)
$$
uniformly on $K$, and the lemma follows since $\psi_{t_n} \to id$ uniformly on $K$ as $n \to \infty$. $\diamond$

\medskip

\begin{lemma} \label{gdot} Let $\{ t_n \}$ be a sequence decreasing monotonically to $0$ such that
$(\phi_{t_n} - id)/t_n$ converges uniformly on compacts of $\dinf$
to a holomorphic function $\chi$. Then $(g_{t_n} - g)/t_n$ converges uniformly on compacts
in a neighbourhood $V$ of $S^1$ to a holomorphic function $\dot{g}$ given (in $V \cap \dinf$) by
$$
\dot{g} = g' \cdot \chi - \chi \circ g
$$
\end{lemma}

\medskip

\noindent{\bf Proof:} By the previous Lemma, the equality
$\psi_{t_n} \circ g = g_{t_n} \circ \psi_{t_n}$ gives
$$
(id - t_n \chi + o(t_n)) \circ g = g_{t_n} \circ (id - t_n \chi + o(t_n))
$$
uniformly on a circle $\{ |z| = R \} \subset \dinf$, so
$$
g - t_n \chi \circ g + o(t_n) = g_{t_n} - g'_{t_n} t_n \chi + o(t_n)
$$
thus
$$
\frac{g_{t_n} - g}{t_n} = g'_{t_n} \cdot \chi - \chi \circ g + o(1) \to g' \cdot \chi - \chi \circ g
$$
uniformly on $\{ |z| = R \}$. Since the maps $g_{t_n}, g$ are circle maps, $(g_{t_n} - g)/t_n$
converges uniformly on the circle $\{ |z| = 1/R \} \subset \D$, hence by the Maximum Principle
we have uniform convergence on the closed annulus $\{ 1/R \leq |z| \leq R \}$.
$\diamond$

\medskip

\begin{theorem} \label{uniqlimit} Let $\{ t_n \}$ be a sequence decreasing monotonically to $0$ such that
$(\phi_{t_n} - id)/t_n$ converges uniformly on compacts of $\dinf$
to a holomorphic function $\chi$. Then $H = \chi/id$ satisfies $\Re H > 0, H(\infty) = 1$,
and is the Herglotz transform $\HH \mu$ on $\dinf$ of $\mu = r^* \mu_{2,g}$, where $r : S^1 \to S^1$
denotes complex conjugation and $\mu_{2,g}$ is the unique $2$-conformal
measure of $g$.
\end{theorem}

\medskip

\noindent{\bf Proof:} Let $\chi_n(z) = (\phi_{t_n}(z) - z)/t_n$ and let $H_n(z) = \chi_n(z)/z$
for $|z| > 1$. Since for $t > 0$, $\Omega_t$ is a proper subdomain of $\dinf$,
by the Schwarz Lemma $|\phi_{t}(z)/z| > 1$ for all $t > 0$. The Moebius map $w = (z-1)/(z+1)$ maps
$\dinf$ conformally to the half-plane $\{ \Re w > 0 \}$, so
$$
\Re\left(\frac{t_n z H_n(z)}{\phi_{t_n}(z) + z}\right) =  \Re\left( \frac{\phi_{t_n}(z) - z}{\phi_{t_n}(z) + z} \right) > 0
$$
and letting $n \to \infty$ it follows that (using $H_n(z) \to H(z), z/(\phi_{t_n}(z) + z) \to 1/2$) that
$\Re H(z) \geq 0$ for all $z$. Moreover $H_n(\infty) = (e^{t_n} - 1)/t_n$,
and $H_n$ converges uniformly to $H$ on compacts in $\dinf$, so $H(\infty) = 1$. If $\Re H(z_0) = 0$ for some
$z_0$ then by the open mapping theorem we must have $H(z) \equiv H(z_0) \in i\R$, contradicting $H(\infty) = 1$, hence
$\Re H(z) > 0$ for all $z$. By the Herglotz Theorem, $H$ is the Herglotz transform on $\dinf$ of a probability
measure $\mu$ on $S^1$. Since $r^2 = id$ it suffices to show that $r^* \mu$ is $2$-conformal for $g$.

\medskip

Let $\mu = f d\lambda + \mu_s$ be the decomposition of $\mu$ into absolutely continuous and singular parts
with respect to Lebesgue measure.

\medskip

By Fatou's Theorem, the radial limits of $\Re H$ exist for a.e. $\xi \in S^1$, and equal $f \circ r$.
From Lemma \ref{gdot} we have $\dot{g} + \chi \circ g = g' \cdot \chi$,
and so
$$
\frac{\dot{g}}{g}(z) + (H \circ g)(z) = \frac{z g'(z)}{g(z)}H(z)
$$
By Lemma \ref{gdot}, the function $\dot{g}$ is holomorphic in a neighbourhood of $S^1$, and, since
$g_{t_n}(\xi)$ converges to $g(\xi)$ along $S^1$ for $\xi$ in $S^1$, we have $\Re (\dot{g}/g)(\xi) = 0$.
Moreover $\xi g'(\xi)/g(\xi) = |g'(\xi)|$ since $g$ is a circle map, so taking real parts in the
equation above as $z$ tends to $\xi \in S^1$ radially gives $(f \circ r \circ g)(\xi) = |g'(\xi)| f(r(\xi))$,
thus
$$
g^*(r^*(f d\lambda)) = g^*((f \circ r)d\lambda) = (f \circ r \circ g)|g'| d\lambda = |g'|^2 (f \circ r) d\lambda
$$
thus the absolutely continuous part $r^*(f d\lambda) = (f \circ r) d\lambda$ of $r^* \mu$
is $2$-conformal.

\medskip

It remains to show that the singular part $\nu := r^* \mu_s$ of $r^* \mu$ is $2$-conformal. Taking imaginary parts in
$(\dot{g}/g)(z) + (H \circ g)(z) = \frac{z g'(z)}{g(z)}H(z)$ as $z$ tends to a point $\xi \in S^1$
radially gives
$$
(\dot{g}/g)(\xi) + i(Q \circ g)(\xi) = |g'(\xi)| iQ(\xi)
$$
The function $\dot{g}/g$ is holomorphic in a neighbourhood of $S^1$, hence bounded on $S^1$,
so $||Q \circ g| - |g'| |Q|| < M$ for some $M > 0$.

\medskip

For $t \gg 1$, let $\nu_t$ denote the measure $\frac{\pi}{2}t 1_{|Q| > t} d\lambda$, so $\nu_t$
converges weakly to $\nu$ as $t \to +\infty$.

\medskip

Given $\xi_0 \in S^1$ and a small $\epsilon > 0$,
let $U \subset S^1$ be a small interval around $\xi_0$ in $S^1$ such
that $(1 - \epsilon) < |g'(\xi)/g'(\xi_0)| < (1 + \epsilon)$ for $\xi$ in $U$.
 For $t \gg 1$ such that $t/(t-M) < 1+\epsilon, t/(t+M) > 1-\epsilon$, we have
\begin{align*}
\nu_t(g(U)) & = \frac{\pi}{2}t \lambda(g(U) \cap \{|Q| > t\}) \\
            & = \frac{\pi}{2}t \lambda(g(\{ \xi \in U, |Q(g(\xi))| > t \})) \\
            & \leq \frac{\pi}{2}t (1 + \epsilon) |g'(\xi_0)| \lambda(\{ \xi \in U, |Q(g(\xi))| > t \}) \\
            & \leq \frac{\pi}{2}t (1 + \epsilon) |g'(\xi_0)| \lambda(\{ \xi \in U, |g'(\xi)||Q(\xi)| > t - M \}) \\
            & \leq t (1 + \epsilon) |g'(\xi_0)| \frac{1}{T} \nu_T(U) \quad \left(\hbox{where } T = \frac{t - M}{(1+\epsilon)|g'(\xi_0)|} \right)\\
            & = (1+\epsilon)^3 |g'(\xi_0)|^2 \nu_T(U) \\
\end{align*}
Similarly we have
$$
\nu_t(g(U)) \geq (1 - \epsilon)^3 |g'(\xi_0)|^2 \nu_{T'}(U)
$$
where $T' = \frac{t+M}{(1 - \epsilon) |g'(\xi_0)|}$.

\medskip

Douady and Yoccoz (\cite{douadyyoccoz}) show that if $\log |g'|$ is of bounded
variation (which is the case if $g$ is $C^2$) then any conformal measure for $g$ has no
atoms. In particular $\nu = r^* \mu_{2,g}$ has no atoms, hence $\nu_t(U) \to \nu(U)$ for
any interval $U$. From the above it follows that given $\epsilon$, there is a $\delta = \delta(\epsilon) > 0$
such that if $U$ is any interval
of length less than $\delta$, and $\xi_0 \in U$,
then
$$
(1 - \epsilon)^3 |g'(\xi_0)|^2 \nu(U) \leq \nu(g(U)) \leq (1 + \epsilon)^3 |g'(\xi_0)|^2 \nu(U)
$$
Given any interval $V$ in $S^1$ and $\epsilon > 0$, let $U_1, \dots, U_n$ be a partition of $V$
into intervals of length less than $\delta(\epsilon)$ centered around points $\xi_1, \dots, \xi_n$.
We may assume that the variation of $|g'|^2$ on each $U_i$ is less than $\epsilon$. Then
\begin{align*}
\nu(g(V)) & = \sum_i \nu(g(U_i)) \\
          & \leq \sum_i (1+\epsilon)^3 |g'(\xi_i)|^2 \nu(U_i) \\
          & \leq (1+\epsilon)^3 \sum_i \int_{U_i} (|g'(\xi)|^2 + \epsilon) d\nu(\xi) \\
          & = (1 + \epsilon)^3 \left(\int_V |g'(\xi)|^2 d\nu(\xi) + \epsilon \right) \\
\end{align*}

Letting $\epsilon$ tend to $0$ gives $\nu(g(V)) \leq \int_V |g'|^2 d\nu$, and similarly
$\nu(g(V)) \geq \int_V |g'|^2 d\nu$. It follows that $\nu = r^* \mu_s$ is $2$-conformal
for $g$, hence $r^*\mu$ is $2$-conformal for $g$ and $\mu = r^* \mu_{2,g}$. $\diamond$

\bigskip

\section{Proofs of main results}

\medskip

We can now prove the main results from the Introduction.

\medskip

\noindent{\bf Proof of Theorem \ref{meastwoconform}:} Given $g \in \dfal$, by uniqueness
of the $2$-conformal measure $\mu_{2,g}$ of $g$, it follows from Theorem \ref{uniqlimit} any normal limit $\chi$ of
the functions $(\phi_t - t)/t, t > 0$ (where $\phi_t : \dinf \to \Omega_t$ is the
normalized Riemann mapping as before) is given by $\chi(z) = z \cdot H(z)$, where $H$ is
the Herglotz transform on $\dinf$ of $r^* \mu_{2,g}$. Since the normal limit is unique,
$$
\frac{\phi_t(z) - z}{t} \to z \cdot H(z)
$$
uniformly on compacts in $\dinf$ as $t \to 0$.

\medskip

For any $t > 0$ and $s > t$, the map $\phi^{-1}_t \circ \phi_s$ is the normalized Riemann
mapping from $\dinf$ to the complement of the unique hull of logarithmic capacity $s - t$
associated to the circle map $g_t$, so it follows from the same argument as above
(applied to $g_t$) that
$$
\frac{\phi^{-1}_t \circ \phi_s(z) - z}{s - t} \to z \cdot H_t(z)
$$
uniformly on compacts in $\dinf$ as $s \to t$, where $H_t$ is the Herglotz transform on $\dinf$
of $r^* \mu_{2, g_t}$. $\diamond$

\bigskip

\noindent{\bf Proof of Theorem \ref{infgenexists}:} Given $g \in \dfal$ and $g_t = \Phi_t(g)$,
since $(\phi_t - id)/t$ converges uniformly on compacts in $\dinf$ to a unique function $\chi$ as $t \to 0$, it follows
from Lemma \ref{gdot} that for any sequence $\{t_n\}$ converging to $0$, the functions
$(g_{t_n} - g)/t_n$ converge uniformly on compacts in a neighbourhood $V$ of $S^1$ to a
holomorphic function $\dot{g}$ on $V$ satisfying (in $V \cap \dinf$)
$$
\dot{g} = g' \cdot \chi - \chi \circ g
$$
So all normal limits of the normal family $\{(g_t - g)/t\}$ coincide, thus
$$
\frac{g_t - g}{t} \to g' \cdot \chi - \chi \circ g
$$
as $t \to 0$, where $\chi(z) = z \cdot H(z)$, $H$ being the Herglotz transform on $\dinf$
of $r^* \mu_{2,g}$. $\diamond$

\bigskip

\noindent{\bf Proof of Theorem \ref{fwduniq}:} Let $g \in \dfal$ and let $\{g_t\}_{0 \leq t < \epsilon} \subset \dfal$
be such that $g_0 = g$ and such that
the right-hand derivatives exist uniformly on a neighbourhood $V$ of $S^1$, and satisfy
$$
\dot{g}_t := \frac{d^+}{ds}_{s = t} g_s(z) = X(g_t)(z) = g'(z) \cdot \chi_t(z) - \chi_t \circ g_t(z)
$$
where $\chi_t(z) = z H_t(z)$, $H_t(z)$ being the Herglotz transform on $\dinf$ of the measure $r^* \mu_{2, g_t}$.
The function $p(z,t) := H_t(z)$ satisfies $\Re p(z,t) > 0, p(\infty, t) = 1$. The map $g_t$
depends continuously on $t$ with respect to the topology of uniform convergence in a neighbourhood of $S^1$, and
hence, being analytic maps, with respect to $C^1$ convergence on $S^1$, so by \cite{douadyyoccoz} the measures
$r^* \mu_{2, g_t}$ depend continuously on $t$ for the weak topology, hence $p(z,t) = H_t(z)$ depends
continuously on $t$ for fixed $z$. It follows (\cite{pomm} \S 6.1) that
there exists a Loewner chain $(\phi_t)_{0 \leq t < \epsilon}$, where the maps $\phi_t$ are
normalized conformal mappings from $\dinf$ onto a decreasing family of simply connected domains $\Omega_t$, such that
$\phi_t(z) = e^t z + O(1)$ near $z = \infty$, and the following right-hand derivatives exist
uniformly on compacts of $\dinf$:
$$
\frac{d^+}{ds}_{|s = t} \phi^{-1}_t \circ \phi_s = \chi_t
$$
Let $(A_t)_{0 \leq t < \epsilon}$ be the increasing family of annular compacts containing
$S^1$ given by the complement in $\hat{\C}$ of $\Omega_t$ and its reflection in $S^1$, so that the
hulls $\HH_t := \hat{\C} - \Omega_t$ are given by $\HH_t = \cD \cup A_t$.

\medskip

Let $\tilde{g}_t = \phi_t \circ g_t \circ \phi^{-1}_t$, then $\tilde{g}_t$ depends
continuously on $t$ with respect to the topology of uniform convergence in a neighbourhood
of $S^1$, hence so does $\tilde{g}'_t$. Let $U = V \cap \dinf$ be a one-sided neighbourhood
of $S^1$. For $s > t$, let $h = s - t$, then as $h \to 0$, we have,
uniformly on any compact in $U$,
\begin{align*}
\tilde{g}'_s & = \tilde{g}'_t + o(h), \\
\phi_s & = \phi_t \circ (id + h \chi_t + o(h)) = \phi_t + h {\phi}'_t \chi_t + o(h),  \\
g_s & = g_t + h \dot{g}_t + o(h) \\
\end{align*}
so $\tilde{g}_s \circ \phi_s = \phi_s \circ g_s$ gives:
\begin{align*}
\tilde{g}_s \circ (\phi_t + h {\phi}'_t \chi_t + o(h)) & = (\phi_t + h {\phi}'_t \chi_t + o(h)) \circ (g_t + h \dot{g}_t + o(h)) \\
\Rightarrow \tilde{g}_s \circ \phi_t + h (\tilde{g}'_s \circ \phi_t)({\phi}'_t \chi_t) + o(h) & = \phi_t \circ g_t + h({\phi}'_t \circ g_t)(\chi_t \circ g_t + \dot{g}_t) + o(h) \\
\Rightarrow \tilde{g}_s \circ \phi_t + h (\tilde{g}'_t \circ \phi_t)({\phi}'_t \chi_t) + o(h) & = \tilde{g}_t \circ \phi_t + h({\phi}'_t \circ g_t)(g'_t \chi_t) + o(h) \\
\Rightarrow (\tilde{g}_s - \tilde{g}_t) \circ \phi_t + h(\tilde{g}_t \circ \phi_t)' \chi_t + o(h) & = h(\phi_t \circ g_t)' \chi_t + o(h) \\
\Rightarrow (\tilde{g}_s - \tilde{g}_t) \circ \phi_t + h(\tilde{g}_t \circ \phi_t)' \chi_t + o(h) & = h(\tilde{g}_t \circ \phi_t)' \chi_t + o(h) \\
\end{align*}
from which it follows that for any $z$ in $U$, the right-hand derivative $\frac{d^+}{dt} \tilde{g}_t(z) = 0$
for all $t$. Since $t \mapsto \tilde{g}_t(z)$ is continuous, we have $\tilde{g}_t(z) = \tilde{g}_0(z) = g(z)$
for all $t$.

\medskip

The map $\tilde{g}_t$, being the conjugate of the circle map $g_t$ by the map $\phi_t$ on $U$, maps the annulus
$\phi_t(U)$ to the annulus $\phi_t(g_t(U))$. Both these annuli have $\partial \HH_t$ as one boundary component,
and $\tilde{g}_t = g$ extends analytically across $S^1$ to be univalent in a neighbourhood of $S^1$ containing
$\HH_t$, hence $g$ leaves $A_t$ invariant. It follows that $A_t$ is the unique annular compact of $g$ such that
the hull $\HH_t = \cD \cup A_t$ has logarithmic capacity $t$, hence $g_t = \Phi_t(g)$. $\diamond$

\medskip

\noindent{\bf Proof of Theorem \ref{germintegral}:} Let $f$ be a germ with rotation number $\alpha$.
Let $(K_t)_{-\infty < t < t_0}$ be the $1$-parameter family of Siegel compacts of $f$ parametrized
by their logarithmic capacities, and let $\phi_t : \dinf \to \hat{\C} - K_t$ be normalized
conformal mappings such that $\phi_t(\infty) = \infty, \phi_t(z) = e^t z + O(1)$ near
$z = \infty$. Fix discs $\D_{r} \subset \D_{r_0}$ with $0 < r < r_0$ such that $f$ maps $\D_r$
univalently into $\D_{r_0}$.
Then for $t \ll -1$, the circle map $g^{f}_t$ (given by conjugating $f$ by $\phi^{-1}_t$) is univalent
on $\phi^{-1}(\D_r - K_t)$, which is an annulus in $\dinf$ with modulus tending to $+\infty$ as $t \to -\infty$.
Therefore the family $(g^{f}_t)$ forms a normal family, and any normal limit of $g^{f}_t$ as $t \to -\infty$
must be a circle map univalent on $\dinf$, hence equal to $R_{\alpha}$, so $g^{f}_t \to R_{\alpha}$ as $t \to -\infty$.

\medskip

If two germs of integral curves $(g^{f_1}_t), (g^{f_2}_t)$ are equal, then there is a $t_0 \in \R$
and neighbourhoods $D_1, D_2$ of the origin such that for $t < t_0$
there is a map $h_t$ univalent on $\hat{\C} - K_t(f_1)$ with $h'_t(\infty) = 1$ conjugating
$f_1$ on $D_1 - K_t(f_1)$ to $f_2$ on $D_2 - K_t(f_2)$ (where $K_t(f_i), i=1,2$ denotes the Siegel
compact of $f_i$ of logarithmic capacity $t$). The maps $(h_t)_{t < t_0}$ form a normal family,
any normal limit of which is univalent on $\C^*$ and has derivative $1$ at $\infty$, hence must be
the identity. Thus $h_t \to id$ as $t \to -\infty$, and $f_1 = f_2$.

\medskip

Finally, given a backward integral curve $(g_t)_{-\infty < t \leq c}$, let $\{t_n\}$ be a sequence in $(-\infty, c]$
decreasing to $-\infty$. Perez-Marco shows in \cite{perezmcirclemaps} that for the circle map $g_{t_0}$,
there exists a germ $f_{t_0}$ with a Siegel compact $K$
such that the fundamental construction of \cite{perezmcirclemaps} applied to the pair $(f_{t_0}, K)$ gives the
circle map $g_{t_0}$.
Conjugating by a scaling if necessary, we may assume that $K$ has logarithmic
capacity $t_0$, so $K = K_{t_0}(f_{t_0})$ and $g_{t_0} = g^{f_{t_0}}_{t_0}$.
By Theorem \ref{fwduniq}, for $t \in [t_0, c]$, $g_{t} = \Phi_{t - t_0}(g_{t_0})$,
so pulling back a Herman compact for $g_t$ to the plane of $g_{t_0}$ and then to the plane
of $f_{t_0}$ gives a Siegel compact $K_{t}(f_{t_0})$ for
$f_{t_0}$ of logarithmic capacity $t$ such that $g_t = g^{f_{t_0}}_t$.

\medskip

Similarly for any $n \geq 0$ we obtain a germ $f_{t_n}$ such that $g_t = g^{f_{t_n}}_t$
for all $t \in [t_n, c]$. So $f_{t_n}$ is given by conjugating $g_c$ by a conformal map
$\phi_n : \dinf \to \hat{\C} - K_{t_n}(f_{t_n})$ with $\phi_n(z) = e^c z + O(z)$
near $z = \infty$. The maps $\{\phi_n\}$ are conformal
mappings on $\dinf$ with fixed derivative at $\infty$ and hence form a normal family.
Fix an annulus $U \subset \dinf$ with boundary components $S^1$ and a Jordan curve $\gamma$
such that $g_c$ is univalent on $U$. The maps $\{f_{t_n}\}$ are normalized
univalent functions on the Jordan domains $D_n$ bounded by $\phi_n(\gamma)$, and hence form a
normal family as well (since the Jordan curves $\phi_n(\gamma)$ are equicontinuous). Any normal limit $f$
then satisfies $g_t = g^{f}_t$ for all $t \leq c$. $\diamond$

\medskip

The fact (\cite{perezmcirclemaps}) that any circle map $g$ arises from a pair $(f, K)$, where $f$ is a
germ with a Siegel compact $K$ also gives:

\medskip

\begin{prop} For any $g \in \dfal$, there exists a backward integral curve $(g_t)_{-\infty < t \leq c}$
with $g_c = g$.
\end{prop}

\medskip

\noindent{\bf Proof:} Given a pair $(f, K)$ which gives rise to the circle map $g$, let $(K_t)_{-\infty < t \leq c}$
be the unique family of Siegel compacts of $f$ parametrized by logarithmic capacity, with $K_c = K$. Applying
the fundamental construction of \cite{perezmcirclemaps} to each pair $(f, K_t)$ gives a backward integral
curve $(g_t)_{-\infty < t \leq c}$ with $g_c = g$. $\diamond$

\bigskip

\section{Linearizable maps and conformal radius of linearization domains}

\medskip

Let $g \in \dfal$ be a circle map which is analytically linearizable. For such a map we
have uniqueness in both forward and backward time for any integral curve with initial condition
$g$:

\medskip

\begin{theorem} Let $(g^{i}_t)_{-\infty < t < c} \subset \dfal, i = 1,2$ be two integral curves of $X$ such that
$g^{1}_{c_0} = g^{2}_{c_0} = g$ for some $c_0 < c$. Then $g^{1}_{t} = g^{2}_{t}$ for all $t < c$.
\end{theorem}

\medskip

\noindent{\bf Proof:} By Theorem \ref{fwduniq} we have $g^{1}_t = g^{2}_t$ for $c_0 \leq t < c$. For $t < c_0$,
we have $g = \Phi_{s}(g^{1}_t) = \Phi_{s}(g^{2}_t)$ where $s = c_0 - t > 0$. Let $U \subset \dinf$ be an invariant
annulus for $g$ with boundary components equal to $S^1$ and a Jordan curve $\gamma \subset \dinf$. Let $\phi^{1}, \phi^{2}$
be the normalized conformal mappings defined on $\dinf$ conjugating $g$ to $g^{1}_t, g^{2}_t$ respectively. For $i= 1,2$,
$\phi^{i}(\gamma)$ is an invariant Jordan curve for $g^{i}_t$ in $\dinf$, and it follows that
$g^{i}_t$ is analytically linearizable. So the Herman compact for $g^{i}_t$ which gives rise to
$g$ (on conjugating by $\phi^i$) is a $g^{i}_t$-invariant annulus $A^i \subset \dinf$ with
boundary components equal to $S^1$ and an invariant Jordan curve $\gamma^i$, such that $\phi^i$ maps $\dinf$
conformally to $\hat{\C} - (\D \cup A^i)$ and conjugates the action of $g$ on $S^1$ to that of $g^{i}_t$ on $\gamma^i$.

\medskip

Let $L = {\phi}^2 \circ ({\phi}^{1})^{-1} : \hat{\C} - (\D \cup A^1) \to \hat{\C} - (\D \cup A^2)$, then $L(\infty) = \infty, L'(\infty) = 1$,
and $L$ conjugates the action of $g^{1}_{t}$ on $\gamma^1$ to that of $g^{2}_{t}$ on $\gamma^2$.
Fix a point $z_1 \in \gamma^1$ and let $z_2 = L(z_1) \in \gamma^2$.
For $i = 1,2$ let ${\eta}^i$ be a conformal map from a round annulus $\{ r_i < |z| < 1 \}$ to $A^i$
mapping $S^1$ to $\gamma^i$ such that ${\eta}^i(z_i) = 1$, then ${\eta}^i$ conjugates the rotation $R_{\alpha}$ to $g^i_t$.

\medskip

Suppose $r_1 \geq r_2$. Then the map $\nu := {\eta}^2 \circ ({\eta}^1)^{(-1)}$ maps $A^1$ into $A_2$,
conjugates the action of $g^{1}_t$ on $\gamma^1$ to that of $g^{2}_{t}$ on $\gamma^2$, and $\nu(z_1) = L(z_1) = z_2$,
hence $\nu = L$ on $\gamma^1$ (since the maps $\nu, L$ differ on $\gamma^1$ by post-composition with a
homeomorphism of $\gamma^2$ commuting with $g^{2}_t$, which must be identity if it fixes a point of $\gamma^2$).

\medskip

It follows that $L$ extends to a univalent map from $\dinf$ into $\dinf$
(by setting $L = \nu$ on $A_1 \cup \gamma^1$). Since $L'(\infty) = 1$, the Schwarz Lemma
implies that $L = id$, hence $g^{1}_t = g^{2}_t$.

\medskip

A similar argument works if $r_2 \geq r_1$. $\diamond$

\medskip

We recall that the {\it conformal radius} $r(D, z_0)$ of a simply connected domain $D$ with a
basepoint $z_0 \in D$ is defined by $r(D, z_0) = h'(0) > 0$, where
$h : \D \to D$ is a conformal map from the unit disk to $D$
satisfying the normalizations $h(0) = z_0, h'(0) > 0$. Note that a disc of radius $R$ centered around
$z_0$ has conformal radius $R$. If $D$ is a simply connected domain in $\hat{\C}$ with basepoint
$z_0 = \infty$ then the conformal radius is defined to be $r = e^{-t}$ where $t$ is the
logarithmic capacity of the hull $\hat{\C} - D$ (so the domain $\{ |z| > R \}$ has
conformal radius $1/R$).

\medskip

Let $f$ be a linearizable germ, and let $(K_t)_{-\infty < t \leq c}$ be the
family of Siegel compacts of $f$ parametrized by logarithmic capacity. Let $D$ be the Siegel
disk (maximal linearization domain) of $f$, and let $h : \D \to D$ be the normalized
conformal mapping satisfying $h(0) = 0, h'(0) > 0$.
For some $t_0 \leq c$, the interiors of the Siegel compacts $K_t$ for $t \leq t_0$ are
linearization domains $D_t \subset D$ for $f$ bounded by analytic Jordan curves ${\gamma}_t$.
Let $r(t) = r(D_t, 0)$ be
the conformal radius of $D_t$, and let $R = r(D, 0)$ be the conformal radius of the Siegel disk. The
normalized conformal mappings of the domains $D_t$ are given by the maps $h_t : \D \to D_t, w \mapsto
h((r(t)/R)w)$. Let $\Omega_t = \hat{\C} - K_t$, and let $\phi_t : \dinf \to \Omega_t$ be the
normalized conformal map satisfying $\phi_t(\infty) = \infty, \phi'_t(\infty) > 0$.

\medskip

Since $\gamma_t = \partial D_t = \partial \Omega_t$ is an analytic Jordan curve the maps
$h_t, \phi_t$ extend analytically across $S^1$, and define an associated 'welding homeomorphism',
which is the analytic circle map $w_t := {h}^{-1}_t \circ \phi_{t|S^1} : S^1 \to S^1$. The
analytic circle map $k_t := w^{-1}_t$ conjugates the rotation $R_{\alpha}$ to the
circle map $g_t = g^f_t$ (arising from the pair $(f, K_t)$).

\medskip

\begin{lemma} For $t < t_0$,

\smallskip

\noindent (i) The conformal radius $r = r(t)$ of the interior $D_t$ of $\gamma_t$ depends
smoothly on the conformal radius $e^{-t}$ of the exterior $\Omega_t$ of $\gamma_t$.

\smallskip

\noindent (ii) The map $k_t$ depends smoothly on $t$.
\end{lemma}

\medskip

\noindent{\bf Proof:} (i) For $r \in (0, r(t_0))$, let $t = t(r) \in (-\infty, t_0)$
be the logarithmic capacity of the hull $K_t = h(\{|w| \leq r/R\})$.
The parametrizations $h_{|\{|w| = r/R\}}$ of the Jordan curves $\gamma_t$
depend smoothly on $r$, hence the boundary values $\phi_{t|S^1}$ of the normalized
conformal mappings $\phi_t$ depend smoothly on $r$ as well (\cite{lanzapreciso}, Thm 3.4).
Since $\phi'_t(\infty) = e^t$ is given in terms of these boundary values by Cauchy's integral formula,
it follows that $t = t(r)$ depends smoothly on $r$. As $r \mapsto t(r)$ is strictly increasing
(by the Schwarz Lemma), the inverse mapping $t \in (-\infty, t_0) \mapsto r = r(t) \in (0, r(t_0))$
is smooth.

\medskip

\noindent (ii) Since the boundary values $\phi_{t|S^1}$ depend smoothly on $r$ and $t = t(r)$ is smooth by (i),
these boundary values depend smoothly on $t$ as well. By \cite{lanzapreciso}, Thm 3.9, the welding maps $w_t$
depend smoothly on $t$, hence so do their inverses $k_t$. $\diamond$

\medskip

\begin{theorem} Let $H_t$ be the Herglotz transform on $\dinf$ of the measure $r^* \mu_{2,g_t}$, and
let $P_t + iQ_t$ be the boundary values on $S^1$ of $H_t$.
Then for $\xi \in S^1$, we have
$$
(P_t \circ k_t)(\xi) = \frac{r'(t)}{r(t)} \cdot |k'_t(\xi)|
$$
and
$$
\frac{\dot{k_t}(\xi)}{k_t(\xi)} + iQ_t(k_t(\xi)) = 0
$$
(where $\dot{k_t}(\xi)$ denotes the derivative with respect to $t$).
\end{theorem}

\medskip

\noindent{\bf Proof:} For $s = t + \epsilon$ with $\epsilon > 0$ small,
and $w \in \dinf$ close to $S^1$, we have
$$
(k^{-1}_t \circ \phi^{-1}_t \circ \phi_s \circ k_s)(w) = \frac{r(s)}{r(t)}(w)
$$
As $\epsilon \to 0$, we have $\phi^{-1}_t \circ \phi_s = id + \epsilon \chi_t + o(\epsilon)$,
where $\chi_t(z) = z H_t(z)$, and $k_s = k_t + \epsilon \dot{k_t} + o(\epsilon), r(s)/r(t) = 1 + \epsilon r'(t)/r(t) + o(\epsilon)$,
hence
\begin{align*}
k^{-1}_t \circ (id + \epsilon \chi_t) \circ (k_t + \epsilon \dot{k_t})(w) + o(\epsilon) & = \left(1 + \epsilon \frac{r'(t)}{r(t)}\right)w + o(\epsilon) \\
\Rightarrow k^{-1}_t \circ (k_t + \epsilon (\dot{k_t} + \chi_t \circ k_t))(w) + o(\epsilon)& = \left(1 + \epsilon \frac{r'(t)}{r(t)}\right)w + o(\epsilon) \\
\Rightarrow \left(id + \epsilon ((k^{-1})'_t \circ k_t)(\dot{k_t} + \chi_t \circ k_t)\right)(w) + o(\epsilon)& = \left(1 + \epsilon \frac{r'(t)}{r(t)}\right)w + o(\epsilon) \\
\end{align*}
thus
$$
\frac{\dot{k_t}(w) + (\chi_t \circ k_t)(w)}{k'_t(w)} = \frac{r'(t)}{r(t)}w
$$
from which we obtain
$$
\frac{\dot{k_t}(w)}{k_t(w)} + (H_t \circ k_t)(w) = \frac{r'(t)}{r(t)} \cdot \frac{k'_t(w)}{k_t(w)} \cdot w
$$
Since $g_t$ is analytically linearizable, the measure $r^* \mu_{2,g_t}$ is absolutely continuous
with respect to Lebesgue measure and has a smooth density. So for any $\xi \in S^1$,
 letting $w$ tend to $\xi$ radially in the equation above and taking real and imaginary parts gives the equalities asserted
in the theorem (using the fact that $\dot{k_t}(\xi)/k_t(\xi)$ is purely imaginary and $\xi k'_t(\xi)/k_t(\xi) = |k'_t(\xi)|$ since
the maps $k_t$ are circle maps). $\diamond$

\medskip

\bibliography{loewconform}

\begin{thebibliography}{BPM15b}

\bibitem[Bis05]{biswas1}
K.~Biswas.
\newblock Smooth combs inside hedgehogs.
\newblock {\em DISCRETE AND CONTINUOUS DYNAMICAL SYSTEMS, vol. 12, 5}, pages
  853--880, 2005.

\bibitem[Bis08]{biswas2}
K.~Biswas.
\newblock Hedgehogs of hausdorff dimension one.
\newblock {\em Ergodic Theory and Dynamical Systems, vol. 28, 6}, pages
  1713--1727, 2008.

\bibitem[Bis15]{biswas3}
K.~Biswas.
\newblock Positive area and inaccessible fixed points for hedgehogs.
\newblock {\em Ergodic Theory and Dynamical Systems,
  http://dx.doi.org/10.1017/etds.2014.143}, pages 1--12, 2015.

\bibitem[BPM13]{bpm3}
K.~Biswas and R.~Perez-Marco.
\newblock Uniformization of higher genus finite type log-riemann surfaces.
\newblock {\em Preprint, http://arxiv.org/pdf/1305.2339.pdf}, 2013.

\bibitem[BPM15a]{bpm1}
K.~Biswas and R.~Perez-Marco.
\newblock Caratheodory convergence of log-riemann surfaces and euler's formula.
\newblock {\em Contemporary Mathematics, Volume 639, DOI:
  http://dx.doi.org/10.1090/conm/639}, pages 197--203, 2015.

\bibitem[BPM15b]{bpm2}
K.~Biswas and R.~Perez-Marco.
\newblock Uniformization of simply connected finite type log-riemann surfaces.
\newblock {\em Contemporary Mathematics, Volume 639, DOI:
  http://dx.doi.org/10.1090/conm/639}, pages 205--216, 2015.

\bibitem[Brj71]{brjuno}
A.~D. Brjuno.
\newblock Analytical form of differential equations.
\newblock {\em Transactions Moscow Math. Soc. 25}, pages 131--288, 1971.

\bibitem[Che11]{cheritat}
A.~Cheritat.
\newblock Relatively compact siegel disks with non-locally connected
  boundaries.
\newblock {\em Mathematische Annalen, {\bf 349} (3)}, pages 529--542, 2011.

\bibitem[dCP03]{lanzapreciso}
M.~Lanza de~Cristoforis and L.~Preciso.
\newblock Differentiability properties of some nonlinear operators associated
  to the conformal welding in schauder spaces.
\newblock {\em Hiroshima Math. J., No. 33}, pages 59--86, 2003.

\bibitem[DY99]{douadyyoccoz}
R.~Douady and J.~C. Yoccoz.
\newblock Nombre de rotation des diffeomorphismes du dercles et mesures
  automorphes.
\newblock {\em Regular and Chaotic Dynamics, Vol. 4, no. 4}, pages 3--24, 1999.

\bibitem[PM93]{perezmens}
R.~Perez-Marco.
\newblock Sur les dynamiques holomorphes non-linearisables et une conjecture de
  v.i. arnold.
\newblock {\em Annales Scientifiques de l'E.N.S 26}, pages 565--644, 1993.

\bibitem[PM94]{perezmtopology}
R.~Perez-Marco.
\newblock Topology of julia sets and hedgehogs.
\newblock {\em preprint, Universit\'e de Paris-Sud}, 1994.

\bibitem[PM95]{perezminvent}
R.~Perez-Marco.
\newblock Uncountable number of symmetries for non-linearisable holomorphic
  dynamics.
\newblock {\em Inventiones Mathematicae 119}, pages 67--127, 1995.

\bibitem[PM96]{perezmdynamics}
R.~Perez-Marco.
\newblock Hedgehogs dynamics.
\newblock {\em preprint, University of California Los Angeles}, 1996.

\bibitem[PM97]{perezmcirclemaps}
R.~Perez-Marco.
\newblock Fixed points and circle maps.
\newblock {\em Acta Mathematica 179:2}, pages 243--294, 1997.

\bibitem[PM00]{perezmsmooth}
R.~Perez-Marco.
\newblock Siegel disks with smooth boundary.
\newblock {\em preprint}, 2000.

\bibitem[Pol96]{poltoratski}
A.~G. Poltoratski.
\newblock On the distribution of boundary values of cauchy integrals.
\newblock {\em Proc. of the AMS, Vol. 124, No. 8}, pages 2455--2463, 1996.

\bibitem[Pom75]{pomm}
C.~Pommerenke.
\newblock Univalent functions.
\newblock {\em Vandenhoeck und Ruprecht}, 1975.

\bibitem[Sie42]{siegel}
C.~L. Siegel.
\newblock Iteration of analytic functions.
\newblock {\em Ann. Math., 43}, pages 807--812, 1942.

\bibitem[Yoc95]{yoccoz}
J.~C. Yoccoz.
\newblock Petits diviseurs en dimension 1.
\newblock {\em Asterisque 231}, 1995.

\end{thebibliography}
\bibliographystyle{alpha}

\medskip

\noindent Ramakrishna Mission Vivekananda University,
Belur Math, WB-711202, India

\noindent email: kingshook@rkmvu.ac.in

\end{document}